\documentclass[11pt,a4paper,twoside]{amsart}

\usepackage{amsfonts, indentfirst}
\usepackage{ifthen,amssymb}
\usepackage[all]{xypic}
\usepackage{young}
\usepackage{graphicx}

\newtheorem{prop}{Proposition}[section]
\newtheorem{theorem}[prop]{Theorem}

\newtheorem{defi}[prop]{Definition}
\newtheorem{coro}[prop]{Corollary}
\newtheorem{remark}[prop]{Remark}

\newenvironment{rem}{\begin{remark}\rm}{\end{remark}}
\newcommand{\cqd}{\hfill$\Box$}

\newcommand{\quot}[0]{\operatorname{Quot}}
\renewcommand{\deg}[0]{{\it deg}}

\newcommand{\dv}[0]{\operatorname{Div}}
\newcommand{\sym}[0]{\operatorname{Sym}}

\newcommand{\gal}[0]{\operatorname{Gal}}

\title[ Codes,  Horn's problem and Gromov-Witten invariants]{Algebraic codes, Horn's problem and Gromov-Witten invariants }

\author{Alberto Besana}
\author{Cristina Mart{\'\i}nez}  

\address{
Fieldaware, 88 Lower Leeson St, Dublin, Ireland}
\email{alberto.besana@fieldaware.com}
  \address{Hamilton Institute, 
Maynooth University, Maynooth, Co. Kildare, Ireland}
 \email{Cristina.MartinezRamirez@nuim.ie}

\begin{document}
\maketitle
\begin{abstract}
We study the Horn problem in the context of algebraic codes on a smooth projective curve defined over a finite field, reducing the problem to the representation theory of the general linear group $GL(n, \mathbb{F}_{q})$. We characterize the coefficients that appear 
in the Kronecker product of symmetric functions in terms of Gromov-Witten invariants of the Hilbert scheme of points in the plane. In addition we classify all the algebraic codes defined over the rational normal curve.
\end{abstract}

\section{Introduction} 
Let denote by $\mathbb{F}_{p}$ the Galois field of $p$ elements. Any other field $\mathbb{F}$ of characteristic $p$ contains a copy of $\mathbb{F}_{p}$. Any $\mathbb{F}_{p^{n}}$ field extension of $\mathbb{F}_{p}$ is a vector space over $\mathbb{F}_{p}$ of dimension $n$ and  an $(n-1)-$dimensional projective space  $PG(n-1,p)$. 
One can consider field extensions $\mathbb{F}_{q}$ of $\mathbb{F}_{p}$ as $q$ varies through powers of the prime $p$.



Let $V$ be an $n+1$ dimensional vector space over the field $\mathbb{F}_{q}$, we denote by $PG(n,q)$ or $\mathbb{P}(V)$ the $n-$dimensional projective space over it and by $\mathbb{P}^{1}$, the projective line. The set of all subspaces of dimension $r$ in $V$ is a Grassmannian and it is denoted by $\mathcal{G}_{r,n}(\mathbb{F}_{q})$
 or by $PG^{r}(n,q)$. 
 The dual of an $r-$space in $PG(n,q)$ is an $(n-r-1)-$space.
 
Consider the $\mathbb{F}_{q}$-rational points of $\mathcal{G}_{r,n}(\mathbb{F}_{q})$ as a projective system, we obtain a $q-$ary linear code, called the Grassmann code, which we denote $G(r,n)$. The lenght $l$ and the dimension $k$ of $G(r,n)$ are given by the $q$ binomial coefficient $l=\left[
\begin{matrix}
 n \\
 r\\
\end{matrix}
\right]_{q}=\frac{(q^{n+1}-1)(q^{n+1}-q)\ldots (q^{n+1}-q^{r})}{(q^{r+1}-1)(q^{r+1}-q)\ldots (q^{r+1}-q^{r})}$, and $k = {n \choose r}$, respectively.

The aim of the present paper is to study the relation between codes constructed from vector bundles and the representation theory of the general linear group $GL(n,\mathbb{F}_{q})$. Namely, consider the right action of the general linear group $GL(n,\mathbb{F}_{q})$ on $\mathcal{G}_{k,n}(\mathbb{F}_{q})$:

\begin{eqnarray}
 \mathcal{G}_{k,n}(\mathbb{F}_{q})\times GL(n,\mathbb{F}_{q}) & \rightarrow & \mathcal{G}_{k,n}(\mathbb{F}_{q}) \\
 (\mathcal{U},A) & \rightarrow & \mathcal{U}A. \nonumber
\end{eqnarray}
Observe that the action is defined independent of the choice of the representation matrix $\mathcal{U}\in \mathbb{F}_{q}^{k\times n}$.

Let $\mathcal{U}\in \mathcal{G}_{k,n}(\mathbb{F}_{q})$ and $G<GL(n,\mathbb{F}_{q})$ be a subgroup, then $C=\{\mathcal{U}A |\,\, A\in G\}$ is an orbit  in $\mathcal{G}_{k,n}(\mathbb{F}_{q})$ of the induced action.

In order to classify all the orbits we need to classify all the conjugacy classes of subgroups of $GL(n,\mathbb{F}_{q})$. In \cite{BM} we studied cyclic coverings of the projective line that correspond to orbits defined by a cyclic subgroup, that is a subgroup in $GL(n,\mathbb{F}_{q})$ containing a cyclic subgroup 
for some prime number. In particular, we showed that any irreducible cyclic plane cover of the projective line can be given by a prime ideal $$(y^{m}-(x-a_{1})^{d_{1}}\ldots (x-a_{n})^{d_{n}})\subset \mathbb{F}_{q}[x,y].$$ This ideal defines an affine curve in $\mathbb{A}^{2}(\mathbb{F}_{q})$ which has singularities, if there are some $d_{k}>1$ for some $1\leq k \leq n$. But there exists an unique smooth projective curve birationally equivalent to this affine curve obtained by homogenization of the polinomial. 
Here we study the connection between ideal sheaves on $\mathbb{F}_{q}[x,y]$ and its numerical invariants together with the combinatorics of partitions of $n$ and the representation theory of the general linear group $GL(\mathbb{F}_{q},n)$. In other words, we want to understand which subspaces are invariant by the action of elements of the general linear group or finite subgroups of $GL(n,\mathbb{F}_{q})$ and how the $GL(n,\mathbb{F}_{q})$ action on the Grassmannian changes the Grassmann code as this action simply permutes basis elements of the Grassmann code.

When one considers as an alphabet a set $\mathcal{P}=\{P_{1},\ldots, P_{N}\}$ of $\mathbb{F}_{q}-$rational points lying on a smooth projective curve defined over a finite field,  algebraic codes are constructed by evaluation of the global sections of a line bundle or a vector bundle on the curve.
Any cyclic cover of $\mathbb{P}^{1}$ which is simply ramified corresponds to an unordered tuple of $n$ points on $\mathbb{P}^{1}$. More generally, in section \ref{sec4} we consider configurations of $n$ points in a $d-$dimensional projective space $PG(d,q)$ which generically lie on a normal rational curve (NRC) and we study the algebraic codes defined on it providing a complete classification in terms of divisors defined over the (NRC), see Theorem \ref{class}. These are the so called Reed-Solomon codes.
Moreover, in the last section as an application of Horn problem we provide a set of generators of the ideal associated to any algebraic code constructed on the NRC over an extension $\mathbb{F}_{q}^{n}$ of $\mathbb{F}_{q}$.




From now $\mathbb{F}_{q}$ will be a field  with $q=p^{n}$ elements and $\mathcal{C}$ a non-singular, projective, irreducible curve defined over  $\mathbb{F}_{q}$ with $q$ elements.

\subsection*{Notation}
For $d$ a positive integer,  $\alpha=(\alpha_{1},\ldots, \alpha_{m})$ is a partition of $d$ into $m$ parts if the $\alpha_{i}$  are positive and decreasing. We will denote as $\mathcal{P}(d)$, the set of all partitions of $d$. We set $l(\alpha)=m$ for the length of $\alpha$, that is the number of cycles in $\alpha$, and $l_{i}$ for the length of $\alpha_{i}$. The notation $(a_{1},\ldots, a_{k})$ stands for a permutation in $S_{d}$ that sends $a_{i}$ to $a_{i+1}$.
A homogeneous symmetric function of degree $n$ over a commutative ring $R$ (with identity) is a formal power series $f(x)=\sum_{\alpha}c_{\alpha}x^{\alpha}$, where $\alpha$ ranges over all weak compositions of $\alpha=(\alpha_{1},\ldots, \alpha_{n})$ of $n$, $c_{\alpha}\in \mathbf{ R}$ and $x^{\alpha}$ stands for the monomial $x^{\alpha_{1}}\cdot x^{\alpha_{2}} \cdots x^{\alpha_{n}}$.
We write $PGL(2,k)=GL(2,k)/k^{*}$, and elements of $PGL(2,k)$ will be represented by equivalence clases of matrices 
$\left(\begin{array}{ll} a &  b \\
c & d \end{array}\right)$, with $ad-bc\neq 0$. 
A $q-$ary constant weight code of length $n$, distance $d$ and weight $w$ will be denoted as an $[n,d,w]_{q}$ code.

\section{Horn problem: an application to convolutional codes}

In this section we present a description of Horn problem in the context of polynomial matrices with polynomial entries associated to torsion modules or dually submodules of a polynomial ring with coefficients in a field. First, we review some known facts about the Quot scheme parametrising quotients of a trivial vector bundle over a generic curve over a field.

Let $\mathcal{O}_{X}$ be the structure sheaf of the curve $X$ defined over a field $k$
and let $K$ be its field of rational functions, considered as a
constant $\mathcal{O}_{X}-$module. Following \cite{BGL}, we define a
divisor of rank $r$ and degree $d$ or $(r,d)$ divisor as a coherent
sub $\mathcal{O}_{X}$-module of $K^{n}= K^{\otimes n}$, having rank
$r$ and degree $d$.
This set can be identified with the set of rational points of an
algebraic variety $\dv^{r,d}_{X/k}$ which may be described as
follows. For any effective ordinary divisor $D$ on $X$, set:
$$\dv^{r,d}_{X/k}(D)=\{E\in \dv^{r,d}_{X,k}; E\subset \mathcal{O}_{X}(D)^{n}\},$$
where $\mathcal{O}_{X}(D)$ is considered as a submodule of $K^{n}$.

The space of all matrix divisors $\mathcal{D}_{k}(r,d)$ of rank $r$ and degree $d$  can be identified with the set of
rational points of $\quot^{m}_{\mathcal{O}_{X}(D)^{n}/X/k}$
parametrizing torsion quotients of $\mathcal{O}_{X}(D)^{n}$ and
having degree $m=r\cdot \deg\,D-d$. It is a smooth projective
irreducible variety. 
Tensoring by
$\mathcal{O}_{X}(-D)$ defines an isomorphism between
$\mathcal{Q}_{r,d}(D)=\quot^{m}_{\mathcal{O}_{X}(D)^{n}/X/k}$ and
$\quot^{m}_{\mathcal{O}^{n}_{X}/X/k}$. 

Since the whole construction is algebraic, it can be performed over any complete valued field $R$ with closed coefficient field $k$ of arbitrary characteristic,
for example, a $p-$adic field or the ring $R=\mathbb{C}\{x\}$ of convergent power series. 
If $f\in R$ is a nonzero divisor, then we define 
the {\it matrix code} $A$ as the code associated to the corresponding torsion module $R/fR$.
The matrix $A$ can be diagonalized by elementary row and column operations with diagonal entries $x^{\alpha_{1}}, x^{\alpha_{2}},\ldots, x^{\alpha_{n}}$, for unique non-negative integers $\alpha_{1} \geq \ldots \geq \alpha_{n}$. 
More precisely, these matrices are in correspondence with endomorphisms of $R^{n},$ with cokernels being torsion modules with at most $n$ generators.  Such a module is isomorphic to a direct sum
$$R/x^{\alpha_{1}}R \oplus R/ x^{\alpha_{2}}R \oplus \ldots \oplus R/x^{\alpha_{n}}R,\ \ \alpha_{1}\geq \ldots \geq \alpha_{n}.$$
The set $(\alpha_{1},\ldots, \alpha_{n})$ of invariant factors of $A$ defines a partition $\alpha$ of size  $d=|\lambda|$.
Reciprocally, when $R=\mathbb{C}\{x\}$ is the ring of convergent power series, any partition $\lambda$ defines a rank one torsion-free sheaf on $\mathbb{C}$ by setting $\mathcal{I}_{\lambda}=(x^{\lambda_{1}},x^{\lambda_{2}}, x^{\lambda_{3}},\ldots, x^{\lambda_{n}})$. In particular, the ideal sheaf corresponding to the identity partition $(1)^{n}$, defines a maximal ideal $\mathcal{I}_{(1)^{n}}=(x,\stackrel{n\  times}{\overbrace{\ldots}},x)$ in $\mathbb{C}[x]$. 
The Horn problem is then equivalent to the following {\bf question:} Which partitions $\alpha, \beta, \gamma$ can be the invariant factors of matrices $A, B$, and $C$ if $C = A \cdot B$? 

In the case of convergent power series, this problem was proposed by I. Gohberg and M.A. Kaashoek. 
Denoting the cokernels of $A, B$ and $C$ by $\mathcal{A}, \mathcal{B}$ and $\mathcal{C}$ respectively, one has a short exact sequence:

$$0\rightarrow \mathcal{A}\rightarrow \mathcal{B} \rightarrow \mathcal{C} \rightarrow 0,$$ i.e.~$\mathcal{B}$ is a submodule of $\mathcal{C}$ with $\mathcal{C}/\mathcal{B}\cong \mathcal{A}$, then such an exact sequence corresponds to matrices $A, B$ and $C$ with $A\cdot B=C$.

If we especialize $C$ to be the identity matrix $I$, 
by the correspondence between partitions and ideal sheaves above, the invariant factors of the identity matrix are defined by the partition $(1)^{n}$, then the question becomes:
Which partitions $\alpha, \beta$ can be the invariant factors of matrices $A$, $B$ if $A\cdot B=I$?
The case of interest for us, will be the case in which $R$ is an $\mathbb{F}_{q}[x]-$module with $q$ a prime power of $p$. 

 Dually, the code can be defined as an $R-$submodule of $R^{n}$, where $R=\mathbb{F}[z]$ is a polynomial ring with coefficients in a field $\mathbb{F}$ and $z$ is a uniformizing parameter in $R$. When $\mathbb{F}$ is a finite field, these are known as convolutional codes which have been very well studied, see for example \cite{CNPS}. 
A full row rank matrix $G(z)\in \mathbb{F}[z]^{k\times n}$ with the property that $$\mathcal{C}={\rm{Im}}_{\mathbb{F}[z]}G(z)=\{f(z)g(z): \, f(z)\in [\mathbb{F}^{k}(z)]\}$$ is called a {\it generator matrix}. The degree $d$ of a convolutional code $\mathcal{C}$ is the maximum of the degrees of the determinants of the $k\times k$ submatrices of one, and hence any, generator matrix of $\mathcal{C}$. The set of convolutional codes of a fixed degree is parametrized by the Grothendieck Quot scheme of degree $d$, rank $n-k$ coherent sheaf quotients of $\mathcal{O}^{n}$ on a curve $X$ defined over $\mathbb{F}$. If the degree is zero, these schemes describe a Grassmann variety and constitute the so called class of block codes of parameters $(n,k)$. 

\subsection{An example with algebraic geometric codes: Reed-Solomon codes}

Let $X$ be a smooth projective curve defined over a finite field $\mathbb{F}_{q}$ with $q$ elements. The classical algebraic-geometric (AG) code due to Goppa is defined by evaluating rational functions associated to a divisor $D$ at a finite set of $\mathbb{F}_{q}-$rational points. From another point of view, we are considering the evaluation of sections of the corresponding line bundle $\mathcal{O}_{X}(D)$ on $X$. Namely, let $\mathcal{P}:=\{P_{1},\ldots, P_{n}\}$ be a configuration of distinct $\mathbb{F}_{q}-$rational points of $X$, the usual algebraic-geometric code is defined to be the image of the evaluation map:

\begin{eqnarray}
\varphi_{D}: L(D)\to \mathbb{F}_{q}^{n} \\
f \mapsto (f(P_{1}),\ldots, f(P_{n})), \nonumber
\end{eqnarray}
where $L(D)$ denotes the vector space of sections associated to the line bundle $\mathcal{O}_{X}$. The parameters of these codes, the length $n$, the dimension $k$ and the minimum distance $d$ are determined by the theorem of Riemann-Roch and it is easy to see that they satisfy the following bound $k+d\geq n+1-g$, where $g$ is the genus of the curve $X$. Using this definition, the notion of AG codes is easily generalised for varieties of higher dimension.

Namely, let $E$ be a vector bundle of rank $r$ on $X$ defined over $\mathbb{F}_{q}$. The Goppa code $C(X,D,G)$ takes as input a divisor $D$ supported on the finite set $\mathcal{P}$ of $\mathbb{F}_{q}-$rational points and a divisor $G$ associated to the vector bundle $E$ and evaluates each section $\sigma \in \mathcal{L}(G)$ in the linear series attached to the divisor $G$:

$$C(X,D,G)=\{(\sigma(P_{i}))_{i=1}^{n}:\, \sigma \in \mathcal{L}(G)\}\subseteq \mathbb{F}^{n}_{q}.$$

Observe that $C(X,\mathcal{P},E)$ is an $\mathbb{F}_{q}-$linear subspace of $\mathbb{F}^{n}_{q^{r}}$ and thus a point of the Grassmannian $\mathcal{G}_{r,n}(\mathbb{F}_{q})$. Moreover, for the same subset of evaluation points and any $r\leq k$, we have $G(r,n)\subseteq G(k,n)\subseteq \mathbb{F}_{q}^{n}$, where $r\leq k$. Further, we get a partial flag of $\mathbb{F}_{q}-$vector spaces $\{0\}=E^{k}\subset E^{k-1}\subset \ldots \subset E^{1}\subset E^{0}=\mathbb{F}^{n}_{q}$ such that $\rm{dim}\, (E^{i-1}/E^{i})=\lambda_{i}$, to which we associate the partition $\lambda=(\lambda_{1},\ldots, \lambda_{r})$ of $n$. In this way, each partition $\lambda$ of $n$ determines a variety $\mathcal{F}_{\lambda}=\mathcal{F}_{\lambda}(\mathbb{F}_{q})$ of partial flags of $\mathbb{F}_{q}-$vector spaces.,

The representation theory of the special linear group $SL(n,\mathbb{F}_{q})$ can be viewed as a form of Gale duality first proven by Goppa in the context of algebraic coding theory.

Let $D$ and $G$ be  effective divisors  supported over a smooth projective curve $X$ defined over $\mathbb{F}_{q}$ such that $\rm{Supp}(G)\cap \rm{Supp}(D)=\emptyset$, then the geometric Goppa code associated with the divisors $D$ and $G$ is defined by $$\mathcal{C}(D,G)=\{(x(P_{1}),\ldots, x(P_{n})), x\in \mathcal{L}(G)\}\subseteq \mathbb{F}_{q}^{n},$$ where $\mathcal{L}(G)$ denotes the linear system associated to the divisor $G$.


\begin{defi}
Let $C_{1}$ and $C_{2}$ be the corresponding codes obtained by evaluating non-constant rational functions $f(x)$ and $g(x)$ with non common roots on $X$ over the support of the divisor $D$. 
We define the {\it quotient code} of $C_{1}$ and $C_{2}$ to be the code associated to the quotient rational function $\varphi=f/g$.
\end{defi}
Since $f$ and $g$  take the value $\infty$, they are defined by non constant polynomials $f(x)$ and $g(x)$ in $\overline{\mathbb{F}_{q}}[x]$. Here $\overline{\mathbb{F}_{q}}$ denotes the algebraic closure of ${\mathbb{F}_{q}}$.
The degree of $\varphi$ is defined to be $\deg\,(\varphi)=\rm{max}\,\{\deg(f), \deg(g)\}$.

As $\varphi$ is a finite morphism, one may associate to each rational point $x\in X(\mathbb{F}_{q})$ a local degree or multiplicty $m_{\varphi}(x)$ defined as:
$$m_{\varphi}(x)=ord_{z=0}\psi(z),$$
where $\psi=\sigma_{2}\circ \varphi \circ \sigma_{1}$,  $y=\varphi(x)$, and $\sigma_{1}, \sigma_{2}\in PGL(2,\mathbb{F}_{q})$ such that $\sigma_{1}(0)=x$ and $\sigma_{2}(y)=0$.

To each non-constant rational function $\varphi$ over $X$, one can associate a matrix $A$ with entries in the ring $\mathbb{F}_{q}[x]$. Namely, let us call $f_{0}:=f(x)$ and call $f_{1}$ the divisor polynomial $g(x)$, and $f_{2}$ the remainder polynomial, then by repeated use of the Euclid's algorithm, we construct a sequence of polynomials $f_{0}, f_{1},\ldots, f_{k}$, and quotients $q_{1},\ldots q_{k}$, $K\leq n$. Then the quotient matrix $A$ is defined to be  the diagonal matrix with entries $q_{1},\ldots, q_{k}$ corresponding to the continued fraction expansion of the rational function~$\varphi$.

Here  we include a SAGE code \cite{SAGE} which implements the algorithm.

{\tt
\begin{samepage}
\begin{verbatim}
def euclid(f, g):
    r = f % g
    q = f // g
    while r.degree() >= 0:
        yield q
        f = g
        g = r
        r = f % g
        q = f // g
\end{verbatim}
\end{samepage}
}





Let $\lambda_{i}$ be the partition of the integer $k$, defining the degree multiplicities of the polynomial $q_{i}$. Then the Horn problem applied to this situation reads:

{\it Which partitions $\alpha, \beta, \gamma$ can be the degree multiplicities of polynomials $q_{A}, q_{B}$ and $q_{C}$ such that the corresponding diagonal matrices $A, B$, and $C$ satisfy
$C = A \cdot B$? }

 Another important family of Goppa codes is obtained considering the rational normal curve $\mathcal{C}^{n}$ defined over $\mathbf{F}_{q}$:
 $$\mathcal{C}^{n}:=\{\mathbb{F}_{q}(1,\alpha, \ldots,\alpha^{n}):\, \alpha\in \mathbb{F}_{q}\cup \{\infty\}\}.$$
The points are distinct elements $\mathbf{F}_{q}$ and $L$ is the vector space of polynomials of degree at most $k-1$ and with coefficients in $\mathbf{F}_{q}$. Such polynomials have at most $k-1$ zeros so nonzero codewords have at least $n-k+1$ non-zeros.  Hence this is a $[n,k,n-k+1]_{q}$ code whenever $k\leq n$. Any codeword $(c_{0},c_{1},\ldots, c_{n-1})$ can be expressed into a $q-$ary $k-$vector with respect to the basis $\{1,\alpha, \ldots, \alpha^{k-1}\}$.These codes are just generalized Reed-Solomon codes  of parameters $[n,k,d]_{q}$ over $\mathbf{F}_{q}$ with parity check polynomial $h(x)=\prod_{i=1}^{q}(x-\alpha^{i})$ where $\alpha$ is a primitive root of $\mathbf{F}_{q}$ such that $\alpha^{k}=\alpha+1$. In other words, the GRS code is an ideal in the ring $\mathbf{F}_{q}[x]/(x^{k}-x-1)$ generated by a polynomial $g(x)$ with roots in the splitting field $\mathbf{F}_{q}^{l}$ of $x^{k}-x-1$, where $k|q^{l}-1$. Since the NRC is a genus 0 curve, it is easy to see that these codes satisfy the Singleton bound $d\geq n-k+1$.
 
Construction of Reed-Solomon codes over $\mathbb{F}_{q}$ only employ elements of $\mathbb{F}_{q}$, hence their lengths are at most $q$. In order to get longer codes, one can make use of elements of an extension of $\mathbb{F}_{q}$, for instance considering subfield subcodes of Reed-Solomon codes. 



As in \cite{BM}, where we considered a variant of the Horn problem in the context of cyclic coverings of the projective line defined over an arbitrary field $k$, the problem is reduced to study the representation theory of the general linear group $GL(n,\mathbb{F}_{q})$.

\section{Representation theory of $GL(n,\mathbb{F}_{q})$}
The multiplication in the finite field $\mathbb{F}_{q^{n}}$ is a bilinear map from $\mathbb{F}_{q^{n}}\times \mathbb{F}_{q^{n}}$ into $\mathbb{F}_{q^{n}}$. Thus it corresponds to a linear map from the tensor product $m: \,\mathbb{F}_{q^{n}}\otimes \mathbb{F}_{q^{n}}\rightarrow \mathbb{F}_{q^{n}}$. The symmetric group $S_{n}$ acts on $\mathbb{F}_{q^{n}}$ via the permutation matrix: 
\begin{eqnarray}
\sigma\cdot v_{i} \ = \  v_{\sigma (i)}, \ \ \ v_{i}\in \mathbb{F}_{q^{n}}. 
\end{eqnarray}
The $d-$Veronese embedding of $\mathbb{P}^{n}(\mathbb{F}_{q})$ maps the line spanned by the vector $v\in \mathbb{F}_{q^{n}}$ to the line spanned by $v^{\otimes d}=v\otimes \ldots \otimes v$. Thus the symmetric group $S_{n}$ acts diagonally on the basis of simple tensors of $\mathbb{F}_{q^{n}}$.

\begin{eqnarray}
 \sigma\cdot (v_{i_{1}}\otimes \ldots \otimes v_{i_{r}})=v_{\sigma(i_{1})}\otimes \ldots \otimes v_{\sigma(i_{r})}.
\end{eqnarray}

For each partition $\lambda=(\lambda_{1},\ldots, \lambda_{k})$ we consider its Young diagram.  The diagram of $\lambda$ is an array of boxes, lined up at the left, with $\lambda_{i}$ boxes in the $i^{th}$ row, with rows arranged from top to botton. For example,
\begin{equation*}
\begin{Young}
 &  &   &   & \cr
 &  & \cr
 &  &  \cr
 \cr
\end{Young}
\end{equation*} is the Young diagram of the partition $\lambda=(5,3,3,1)$ with $l(\lambda)=4$ and 
$|\lambda|=12$. We define the Schur projection:
$$c_{\lambda}: \, \bigotimes^{d}\, \mathbb{F}_{q^{n}}\rightarrow \bigotimes^{d} \mathbb{F}_{q^{n}}.$$

Let $S_{n}$ be the symmetric group of permutations over $d$ elements. Any permutation $\sigma \in S_{n}$ acts on a given Young diagram by permuting the boxes. Let $R_{\lambda}\subseteq S_{n}$ be the subgroup of permutations preserving each row. Let $C_{\lambda}\subseteq S_{n}$ be the subgroup of permutations preserving each column, let $c_{\lambda}= \sum_{\sigma \in R_{\lambda}} \sum_{\tau \in C_{\lambda}} \epsilon (\tau)\sigma \tau$. 


 

The image of $c_{\lambda}$ is an irreducible $GL(n,\mathbb{F}_{q})-$module, which is nonzero iff the number of rows is less or equal than 
dim$V_{\lambda}$. All irreducible $GL(n,\mathbb{F}_{q})-$mo\-du\-les can be obtained in this way. Every $GL(n,\mathbb{F}_{q})-$module is a sum of irreducible ones.

In terms of irreducible representations of
$GL(n,\mathbb{F}_{q})$, a partition $\eta$ corresponds to a finite
irreducible representation that we denote as $V(\eta)$. Since
$GL(n,\mathbb{F}_{q})$ is reductive, any finite dimensional
representation decomposes into a direct sum of irreducible
representations, and the structure constant $c^{\eta}_{\lambda, \mu}$
is the number of times that a given irreducible representation
$V(\eta)$ appears in an irreducible decomposition of
$V(\lambda)\otimes V(\mu)$. These are known as Littlewood-Richardson
coefficients, since they were the first to give a combinatorial
formula encoding these numbers (see \cite{Fu}). 
In terms of the Hopf algebra $\Lambda$ of Schur functions, let $s_{\lambda}$ be the Schur function indexed by the partition $\lambda$, we have  $s_{\lambda}\cdot s_{\mu}= \sum_{\nu} c_{\lambda\mu}^{\nu}s_{\nu}$ for the product and we get the coefficients $k_{\lambda \mu}^{\eta}$ as the structure constants of the dual Hopf algebra $\Lambda^{*}$. These are known as Kronecker coefficients, (see \cite{Ma} and \cite{SLL}). Recall that the Schur function $s_{\lambda}$ attached to the partition $\lambda=(\lambda_{1},\ldots, \lambda_{n})$ of lenght less or equal than $n$ is defined by the quotient:
$$s_{\lambda}(x_{1},\ldots, x_{n})=\frac{det(x_{i}^{\lambda_{i}+n-j})_{1\leq i,j, \leq n}}{det(x^{n-j}_{i})_{1\leq i,j\leq 0}}.$$ It is a homogeneous polynomial of degree $|\lambda|$ in $x_{1},\ldots, x_{n}$. It easily seen that $s_{\lambda}(x_{1},\ldots, x_{n},0)=s_{\lambda}(x_{1},\ldots, x_{n})$. Moreover we can define the Schur function $s_{\lambda}$ as the unique symmetric function with this property for all $n\geq l(\lambda)$. It is well known that the Schur functions constitute a basis for the ring $\Lambda$ of symmetric functions. 
In addition, there are at least other three well known bases for the ring $\Lambda$ of symmetric functions. The basis $e_{k}$ of $k-$elementary symmetric functions,  the $h_{k}$ complete homogeneous symmetric functions of degree $k$ and the power sums $p_{k}=z^{k}_{1}+z^{k}_{2}+\ldots$. 
This has been applied in Reed-Solomon coding, that is, for AG codes defined on the projective line $\mathbb{P}^{1}$, as a way to encode information words. 
Namely, for each codeword $a=(a_{0},a_{1},\ldots, a_{n})$, $a_{i}\in \mathbb{F}_{q}$, let us define $a_{n+1}=\sum_{i=1}^{n}a_{i}\in \mathbb{F}_{q}$ which is nothing but the first elementary symmetric function $e_{1}$. If we consider the variables $x_{1},\ldots, x_{r}$ as a fixed list of nonzero elements in $\mathbb{F}_{q}$, then the information word $a$ can be encoded into the codeword $d=(d_{1},\ldots, d_{r})$, where $d_{i}=\sum_{j=1}^{n}a_{j}x_{i}^{j}$. The secret is $a_{0}=-\sum_{i=1}^{r}d_{i}$, while the pieces of the secret are the $d_{i}'s$. 

\subsection{Relation between Littlewood-Richardson coefficients and Kronecker coefficients}
One can stack Littlewood-Richardson coefficients $c^{\nu}_{\lambda \mu}$ in a 3D matrix or 3-dimensional matrix. Intuitively a 3D matrix  is a stacking of boxes in the corner of a room. The elements of the principal diagonal are called rectangular coefficients and are indexed by triples $(\lambda, \mu, \nu)=((i^{n}), (i^{n}),(i^{n}))$ of partitions $(i^{n})$ with all their parts equal to the same integer $1\leq i \leq n$.
\begin{figure}[htb]
\centering
\includegraphics[scale=0.3]{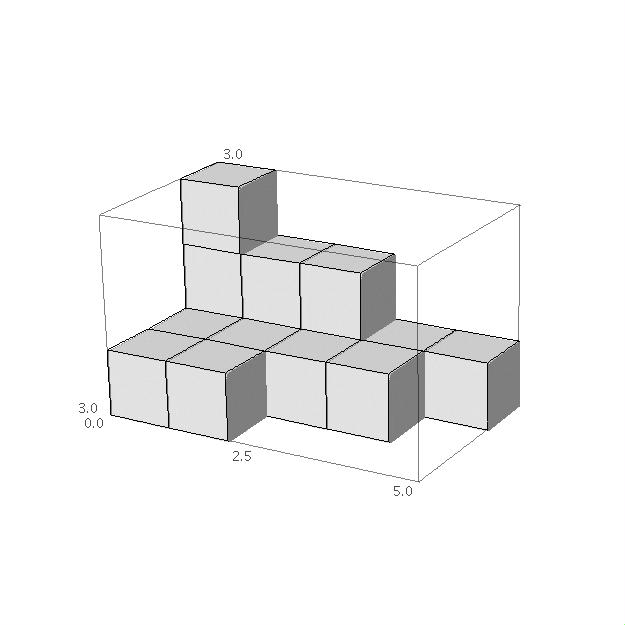}
\end{figure}

Consider $\mathcal{B}$ and $\mathcal{C}$, two 3D matrices, then we define the product matrix $\mathcal{B}\cdot \mathcal{C}$ as the 3D matrix $$\mathcal{B}\cdot \mathcal{C}= \prod_{\nu \in \mathcal{P}(n), B^{\nu},C^{\nu}\in M_{p(n)\times p(n)}(\mathbb{Q})}B^{\nu}\cdot C^{\nu}. $$
Namely, for each index $\nu$ fixed, $\lambda$ and $\mu$ run over all partitions $\mathcal{P}(n)$ of $n$. Thus the coefficients $\left(c^{\nu}_{\lambda, \mu}\right)_{\lambda, \mu \in \mathcal{P}(n)}$ are encoded in a matrix of order $p(n)\times p(n)$, where $p(n)$ denotes the number of unrestricted partitions of $n$, that is, the number of ways of writing the integer $n$ as a sum of positive integers without regard to order. Thus the product matrix $B^{\nu}\cdot C^{\nu}$ is the standard product of square matrices in $M_{p(n)\times p(n)}(\mathbb{Q})$. In particular, the property of associativity follows easily from the associativity in the space $M_{p(n)\times p(n)}(\mathbb{Q})$.


\begin{prop} Let $\mathcal{C}$ be the 3D matrix whose entries are the Littlewood-Richardson coefficients, and $\mathcal{K}$  the 3D matrix of Kronecker coefficients. Then the  matrices are inverse one to each other. 

\end{prop}

{\it Proof.} Since $c^{\nu}_{\lambda \mu}$ and $k^{\nu}_{\lambda \mu}$ correspond to the structure constants of the Hopf algebra of Schur functions and its dual one respectively (see Proposition 4.3 of \cite{BM}), and the Hopf algebra of Schur functions is self-dual (see \cite{SLL}), one gets that the product matrix $\mathcal{C}\cdot \mathcal{K}$ is the identity 3D  matrix $\mathcal{I}$, that is, the matrix whose rectangular coefficients are identically 1.\ 
Thus both matrices are inverse one to each other.
\cqd


\subsection{The polytope of triples $(\lambda, \mu, \eta)$ for which $c^{\eta}_{\lambda, \mu}$ is positive}
The convex hull in $\mathbb{R}^{3}$ of all triples $(\lambda, \mu, \nu)$ with $c^{\nu}_{\lambda, \mu}>0$ is the Newton polytope of $f(x,y,z)=\sum _{\lambda, \mu, \nu} c^{\nu}_{\lambda, \mu}x^{\lambda} y^{\mu} z^{\nu}  \in \mathbb{C}[x,y,z]$. Here $x^{\lambda}$ denotes the monomial $x^{\lambda_{1}}\cdots x^{\lambda_{n}}$ of partition degree $\lambda$. In particular, when $\lambda=(1^{r})$, we have $x^{(1)^{r}}=e_{r}=\sum_{i_{1}<\ldots <i_{r}}\, x_{i_{1}}\ldots x_{i_{r}}$, the $r-$th elementary symmetric function. At the other extreme, when $\lambda=(r)$ we have $x^{(r)}=p_{r}=\sum\, x_{i}^{r}$, the $r-$power sum. As we have seen in the previous the section, it is clear that every symmetric function $f\in \Lambda$ is uniquely expressible as a finite linear combination of the $(x^{\lambda})_{\lambda\in \mathcal{P}}$. Moreover, the following Theorem shows that  $f$ is the  the generating series for the Gromov-Witten invariant $N_{d,g}(\lambda, \mu, \nu)$ counting irreducible plane curves of given degree and genus $g$ passing through a generic configuration of $3d-1+g$ points on $\mathbb{P}^{2}(\mathbb{C})$ with ramification type at $0, \infty$ and 1 described by the partitions $\lambda, \mu$ and $\nu$ and simple ramification over other specified points with $|\lambda|+|\mu|+|\nu|=d,$ and these have been computed by Fomin and Mikhalkin in \cite{FM}.

\begin{theorem} 
The power series 
$f(x,y,z)=\sum _{\lambda, \mu, \nu} c^{\nu}_{\lambda, \mu}x^{\lambda} y^{\mu} z^{\nu}  \in \mathbb{C}[x,y,z],$ is the generating series for the 
Gromov-Witten invariants $N_{d,g}(\lambda, \mu, \nu)$, counting irreducible plane curves of given degree $d$ and genus $g$ passing through a generic configuration of $3d-1+g$ points on $\mathbb{P}^{2}(\mathbb{C})$ with ramification type at $0, \infty$ and 1 described by the partitions $\lambda, \mu$ and $\nu$ and simple ramification over other specified points with $|\lambda|+|\mu|+|\nu|=d$.%
\end{theorem}

{\it Proof.} Whenever the coefficient $c^{\nu}_{\lambda, \mu}>0$ is positive consider the corresponding ideal sheaves $\mathcal{I}_{\lambda}$, $\mathcal{I}_{\mu}$ and $\mathcal{I}_{\nu}$ in $\mathbb{C}[x,y,z]$ associated to the partitions $\lambda, \mu$ and $\nu$ respectively. Each ideal sheaf determines a curve in $\mathbb{C}[x,y]$ via homogenization of the corresponding monomial ideals.
Thus each coefficient represents the number of ideal sheaves on $\mathbb{C}^{3}$ of colength $n$ and degree $d$ equal to the size of the partition, that is the corresponding 
3-point  Gromov-Witten invariant $\langle \lambda, \mu, \nu \rangle _{0,3,d}$
of the Hilbert scheme ${\rm{Hilb}}_{n}$ of $n=2d-1+|\nu|+|\mu|+|\lambda|+g$ distinct points in the plane, or the relative Gromov-Witten invariant $N_{d,g}(\lambda, \mu, \nu)$ counting irreducible plane curves of given degree $d$ and genus $g$ passing through a generic configuration of $3d-1+g$  points on $\mathbb{P}^{2}(\mathbb{C})$  with ramification type at $0, \infty$ and 1 respectively, described by the partitions $\lambda, \mu$ and $\nu$ of $n$, (see section 4 of \cite{BM}). 
\cqd

\begin{rem} 
The Euler characteristic of each ideal sheaf is fixed and coincides with the Euler characteristic $\chi$ of the polyhedra described in $\mathbb{R}^{3}$ by the convex hull  of all triples $(\lambda, \mu, \nu)$ with $c^{\nu}_{\lambda, \mu}>0$, that is,  the Newton polytope of $f(x,y,z)=\sum _{\lambda, \mu, \nu} c^{\nu}_{\lambda, \mu}x^{\lambda} y^{\mu} z^{\nu}  \in \mathbb{R}[x,y,z].$
Thus each coefficient represents the number of ideal sheaves on $\mathbb{C}^{3}$ of fixed Euler characteristic $\chi=n$ and degree $d$ equal to the size of the partition, that is the corresponding 
Donaldson-Thomas invariant of the blow-up of the plane $\mathbb{P}^{1}\times (\mathbb{C}^{2})$ with discrete invariants  $\chi=n$ and degree $d$.

\end{rem}

\begin{rem}The Hilbert scheme ${\rm{Hilb}}_{n}$ of $n$ points in the plane $\mathbb{C}^{2}$ pa\-ra\-me\-tri\-zing ideals $\mathcal{J}\subset \mathbb{C}[x,y]$ of colength $n$ 
contains an open dense set in the Zariski topology para\-me\-tri\-zing ideals associated to configurations of $n$ distinct points. Moreover there is an isomorphism ${\rm{Hilb}}_{n}\cong (\mathbb{C}^{2})^{n}/S_{n}$. In particular, as we showed in \cite{BM}, any conjugacy class in the symmetric group $S_{n}$ determines a divisor class in the $T-$equivariant cohomology $H^{4n}_{T}({\rm{Hilb}}_{n},\mathbb{Q})$, for the standard action of the torus $T=(\mathbb{C}^{*})^{2}$ on $\mathbb{C}^{2}$.  The $T-$equivariant cohomology of ${\rm{Hilb}}_{n}$ has a canonical Nakajima basis indexed by $\mathcal{P}(n)$. The map $\lambda \rightarrow \mathcal{J}_{\lambda}$ is a bijection between the set of partitions $\mathcal{P}(n)$ and the set of $T-$fixed points ${\rm{Hilb}}^{T}_{n}\subset {\rm{Hilb}}_{n}$.
\end{rem}
Denote the series $\langle \lambda, \mu, \nu \rangle^{{\rm{Hilb}}_{n}}$ of 3-point invariants by a sum over curve degrees:
$$\langle \lambda, \mu, \nu \rangle^{{\rm{Hilb}}_{n}}=\sum_{d\geq 0}q^{d} \langle \lambda, \mu, \nu \rangle^{{\rm{Hilb}}_{n}}_{0,3,d}.$$
\begin{coro} Let $H$ be the divisor class in the Nakajima basis corresponding to the tautological rank $n$ bundle $\mathcal{O}/\mathcal{J}\rightarrow {\rm{Hilb}}_{n}$ with fiber $\mathbb{C}[x,y]/\mathcal{J}$ over $\mathcal{J}\in {\rm{Hilb}}_{n}$ and $\nu$ the corresponding partition.
Then we can recover inductively in the degree $d$, all the Littlewood-Richardson coefficients $(c^{\nu}_{\lambda,\mu})_{\lambda, \mu \in \mathcal{P}(n)}$.
\end{coro}

{\it Proof.}
The non-negative degree of a curve class $\beta \in H_{2}({\rm{Hilb}}_{n},\mathbb{Z})$ is defined by $d=\int_{\beta}H$. Then via the indentification of $c^{\nu}_{\lambda,\mu}$ with the 3-point Gromov-Witten invariant $\langle \lambda,H, \mu\rangle^{{\rm{Hilb}}_{n}}_{0,3,d}$ where $[\lambda],[\mu]$ are the corresponding classes in $H^{4n}_{T}({\rm{Hilb}}_{n},\mathbb{Q})$ associated to the partitions $\lambda$ and $\mu$ in $\mathcal{P}(n)$, we proceed by induction on the degree $d$ as in section 3.6 of \cite{OP}. \cqd

\begin{rem} If we choose the partition $\nu$ to be the empty partition $\emptyset$, we recover the relative Gromov-Witten invariants $N_{d,g}(\lambda, \mu)$ studied by Fomin and Mikhalkin in \cite{FM}, and by Caporaso and Harris in \cite{CH}. 
\end{rem}





\section{Configurations of points over a normal rational curve}\label{sec4}
In this section, we study codes defined from a linear series attached to a divisor on the normal rational curve (NRC) or equivalently Goppa codes on $\mathbb{P}^{1}$ and hence generalised-Reed Solomon codes. 
Assume $V$ is a  vector space of dimension $n+1$ over a field $k$ equipped with a linear action, that is, $G$ acts via a representation $G \rightarrow GL(V )$. We denote by $S^{d}V$ the $d-$th symmetric power of $V$. 


Consider the $d-$Veronese embedding of $\mathbb{P}^{n}$ 
\begin{eqnarray}
 \mathbb{P}V^{*}\rightarrow \mathbb{P}S^{d}V^{*}  \\
v \mapsto v^{d}, \nonumber
\end{eqnarray}
 mapping the line spanned by $v\in V^{*}$ to the line spanned by $v^{d}\in S^{d}V^{*}$. In coordinates, if we choose bases $\{\alpha, \beta\}$ for $V$ and $\{[\frac{n!}{k!(n-k)!}]\alpha^{k }\beta^{d-k}\}$ for $S^{d}V^{*}$ and expanding out $(x\alpha+y\beta)^{d}$, we see that in coordinates this map may be given as

$$[x,y]\rightarrow [x^{d},x^{d-1}y, x^{d-2}y^{2}, \ldots, xy^{d-1},y^{d}].$$ 


In particular, the homogeneous coordinate ring for the natural projective embbeding of the geometric invariant theory (GIT) quotient $(\mathbb{P}^{d})^{n}// SL_{d+1}$ is the ring of invariants for $n$ ordered points in the projective space up to projectivity, i.e, if one considers the function field $k(x_{1},\ldots, x_{d})$ of the projective space $(\mathbb{P})^{d}$, the ring of invariants is defined by:
$$\{f\in k(x_{1},\ldots, x_{d})\,|\, \forall \sigma \in SL_{d+1}, \sigma\cdot f=f \}. $$

Generators for this ring are given by tableau functions, which appear in many areas of mathematics, particularly representation theory and Schubert calculus. Consider the hypersimplex:
$$\triangle(d+1,n)=\{(c_{1},\ldots,c_{n})\in \mathbb{Q}^{n}| 0\leq c_{i}\leq 1,  \sum c_{i}=d+1\},$$ for any $1\leq d\leq n-3$ and choose a linearization $c\in \triangle (d+1,n)$, there is a morphism $$\varphi: \bar{M}_{0,n}\rightarrow (\mathbb{P}^{d})^{n} //_{c}SL_{d+1},$$ sending a configuration of distinct points on $\mathbb{P}^{1}$ to the corresponding configuration under the $d^{th}$ Veronese map.

The symmetric power $\sym^{n}C_{d}$ of the curve $C_{d}$ is the quotient of  the configuration space $\mathcal{C}^{n}_{d}$ of $n$ unordered tuples of points on the rational normal curve $C_{d}$ by the symmetric group $S_{n}$. 
Furthermore, we can identify the set of effective divisors of degree $d$ on $C_{d}$ with the set of $k-$rational points of the symmetric power ${\rm{Sym}}^{n}C$, that is, ${\rm{Sym}}^{n}C$ represents the functor of families of effective divisors of degree $n$ on $C$.

\subsection*{Why codes on the normal rational curve?}
By definition, the rational normal curve $C_{d}$ is the image by the $d-$Veronese embedding of $\mathbb{P}V^{*}=\mathbb{P}^{1}$ where $V$ is a 2-dimensional vector space, therefore it is isomorphic to any curve of genus 0. The action of $PGL(2,k)$ on $\mathbb{P}^{d}$ preserves the rational normal curve $C_{d}$. Conversely, any automorphism of $\mathbb{P}^{d}$ fixing $C_{d}$ pointwise is the identity. It follows that the group of automorphisms of $\mathbb{P}^{d}$ that preserves $C_{d}$ is precisely $PGL(2,k)$. These codes are just generalized RS codes and they come with efficient decoding algorithms once we choose a metric consistent with channel errors and search of a set of vectors with given metric properties as a correcting code. In particular, these codes are consistent with the Hamming metric (\cite{BS1}, \cite{BS2}). Recall that given two vectors of length $n$, say $U$ and $V$, the Hamming distance $d_{H}(U,V)$ between $U$ and $V$ is the number of coordinates in which they differ. 

Given a $[n,k]$ RS code $C$ of length $n$ and dimension $k$, we call $d$ the minimum (Hamming distance) which satisfies the Singleton bound $n-k+1$. We shall identify the code with the set of its codewords. A codeword of $C$ is viewed as a polynomial $c_{0}+c_{1}z+\ldots+c_{n-1}z^{n-1}$ in the $\mathbb{F}-$vector space $\mathbb{F}[z],$ where $\mathbb{F}$ is a finite field. 
 In the communication process, when a codeword is transmitted, it can be affected by errors and erasures. An error occurs when one codeword component is changed into another field element and an erasure occurs when the received component has an unknown value. The problem of minimum distance decoding is to find, for any given vector $r$, the set $C_{r}$ of all codewords $c\in C$ at minimum distance from $r$. If $C_{r}$ contains just one element $c$, then the sent codeword coincides with the received codeword and no decoding is needed. The codewords of minimum weight are the points lying in the intersection of any line and the curve. K. Lee and M.E. O'Sullivan in \cite{LS} describe a list decoding algorithm consisting of two steps: the interpolation step and the root-finding step. %
 Starting with a set of generators of the module induced from the ideal for the $n$ points $\{P_{1},\ldots, P_{n}\}$, they convert the generators to a Gr\"obner basis of the module in which the minimal polynomial is found. This results in an efficient algorithm solving the interpolation problem. In section 6, we provide Horn's algorithm to compute sets of indices which are admissible for the Horn problem. As a result, we provide a set of generators for the algebraic code induced on the NRC.




\begin{prop} \label{div} If we consider the set of orbits of $\mathcal{C}^{n}_{d}$ by the action of finite subgroups of the symmetric group $S_{n}$, we get all possible divisor classes in the group ${\rm{Div}}^{n}(C_{d})$ of degree $n$ divisors on $C_{d}$.
\end{prop}
{\it Proof.}
Since the symmetric group $S_{n}$ is generated by 3 elements, a reflection of order 2, a symmetry of order 3 and a rotation of order $n$, we get all the divisor classes by quotienting the configuration space $\mathcal{C}^{n}_{d}$ of $n$ points on the rational normal curve, by the cyclic group generated by the rotation, or one of the triangle groups, the dihedral group $D_{n}$, the alternated groups $A_{4}$, $A_{5}$ or the symmetric group $S_{4}$. \cqd


\section{Notion of collinearity on the rational normal curve}
A permutation matrix $\sigma \in GL(n,\mathbb{F}_{q})$ acts on the Grassmannian by multiplication on the right of the corresponding representation matrix. In particular, we are interested in understanding the orbits by the action of any permutation matrix of $GL(n,\mathbb{F}_{q})$ and moreover of any subgroup $G$ contained in $GL(n,\mathbb{F}_{q})$. Further, it is possible to count the orbits of the action in several cases and these correspond to sets of points satisfying certain geometrical conditions.
\begin{defi}
An incidence structure  $\mathcal{S}$  on $V$ is a triple $(\mathcal{P},\mathcal{B},I)$, where $\mathcal{P}$ is a set whose elements are smooth, reduced points in $V$, $\mathcal{B}$ is a set whose elements are subsets of points called blocks (or lines in several specific cases) endowed with a relation of collinearity, and an incidence relation $I\subset \mathcal{P}\times \mathcal{B}$. 
If $(P,L)\in I$, then we say that $P$ is incident with $L$ or $L$ is incident with $P$, or $P$ lies in $L$ or $L$ contains $P$. 
\end{defi}

When the collinearity relation is a symmetric ternary relation defined on triples $(p,q,r)\in  \mathcal{P}\times \mathcal{P}\times \mathcal{P}$ by the geometric condition $(p,q,r) \in \mathcal{B}$ if either $p+q+r$ is the full intersection cycle of $C_{d}$ with a $k-$line $l\subset \mathbb{P}^{n}(k)$
 with the right multiplicities, or else if there exists a $k-$line $l\subset V$ such that, $p,q,r \in l$, then the triple $(p,q,r)$  is called a plane section.

\begin{enumerate}
\item  For any $(p,q)\in \mathcal{P}^{2}(V^{*})$, there exists an $r\in \mathcal{P}(S^{d}V^{*})$ such that $(p,q,r)\in l$. The triple $(p,q,r)$ is strictly collinear if $r$ is unique with this property, and $p, q, r$ are pairwise distinct. The subset of strictly collinear triples is a symmetric ternary relation. When $k$ is a field algebraically closed of characteristic 0, then $r$ is unique with this property, and we recover the euclidean axioms.
\item Assume that $p\neq q$ and that there are two distinct $r_{1},r_{2}\in \mathcal{P}$ with $(p,q,r_{1})\in \mathcal{B}$ and $(p,q,r_{2})\in \mathcal{B}$. Denote by $l=l(p,q)$ the set of all such $r's$, then $l^{3}\in \mathcal{B}$, that is any triple $(r_{1},r_{2},r_{3})$ of points in $l$ is collinear. Such sets $l$ are called lines in $\mathcal{B}$.

 \end{enumerate}
 
 If $V$ is a 3-dimensional vector space defined over the finite field $\mathbb{F}_{p}$, then 
the projective plane $\mathbb{P}^{2}(\mathbb{F}_{p})$  on $V$ is defined by the incidence structure $PG(2,p)=(\mathcal{P}(V),\mathcal{L}(V),I)$.

\begin{defi}
\begin{enumerate}
 
 \item A $(k;r)-$arc $\mathcal{K}$ in $PG(2,p)$ is a set of $k-$points such that some $r$, but not $r+1$ of them are collinear. In other words, some line of the plane meets $\mathcal{K}$ in $r$ points and no more than $r-$points. A $(k;r)-$arc is complete if there is no $(k+1;r)$ arc containing it.

 \item A $k-$arc is a 
 set of $k$ points, such that, every subset of $s$ points with $s\leq n$ points is linearly independent.
 \end{enumerate}
 
\end{defi}

Let $q$ denote some power of the  prime $p$ and $PG(n,p)$ be the $n-$dimensional projective space $(\mathbb{F}_{p})^{n+1}\cong \mathbb{F}_{q}$, where $n\geq 2$.
The normal rational curve $C$ is defined as:
$$\mathcal{V}^{n}_{1}:=\Big\{\mathbb{F}_{q} (1,x,x^{2},\ldots,x^{n})| \ x\in \mathbb{F}_{q}\bigcup \{\infty\}\Big\}.$$ 




If $q\geq n+2$, the NRC is an example of a $(q+1)-$arc. It contains $q+1$ rational points, and every set of $n+1$ points are linearly independent.
For each $a\in (\mathbb{F}_{p})^{n+1}$, the mapping:
$$\mathbb{F}_{p}(x_{0},\ldots, x_{n})\rightarrow \mathbb{F}_{p}(a^{0}x_{0},\ldots, a^{n}x_{n}),$$ describes an automorphic collineation of the NRC. 

All invariant subspaces form a lattice with the operations of "join" and "meet".

For $j\in \mathbb{N}$, let
$\Omega(j)=\{m \in \mathbb{N}| 0 \leq m\leq n, {m \choose j} \neq 0 \, mod \,p\}$.
Given $J\subset \{0,1,\ldots, n\}$, put
$\Omega(J)=\bigcup_{j\in J}\Omega(j), \, \Psi(J):=\bigcup_{j\in J}\{j,n-j\}$.

Both $\Omega$ and $\Psi$ are closure operators on $\{0,1,\ldots,n\}$. 
Likewise the projective collineation
$\mathbb{F}_{p}(x_{0},x_{1},\ldots, x_{n})\rightarrow \mathbb{F}_{p}(x_{n},x_{n-1},\ldots, x_{0})$ leaves the NRC invariant whence $\Lambda$ has to be closed with respect to $\Psi$. 
Any algebraic-geometric code constructed by evaluation of a function over the NRC with values in  $\mathbb{F}_{q}$ is a generalised Reed-Solomon code of length at most $q$. In order to get longer codes, one needs to use elements from any finite extension $\mathbb{F}_{q}^{r}$ of $\mathbb{F}_{q}$.

\begin{prop} \label{codeNRC}
Each subspace invariant under collineation of the NRC, is indexed by a partition in $\mathcal{P}(t)$. If the ground field $k$ is sufficiently large, then every subspace which is invariant under all collineations of the NRC, is spanned by base points $kc_{\lambda}$, where $\lambda \in \mathcal{P}(t)$. 
\end{prop}

{\it Proof.}
Let $$E^{t}_{n}:=\{(e_{0},e_{1},\ldots, e_{n})\in \mathbb{N}^{n+1}|\, e_{0}+e_{1}+\ldots+e_{n}=t\},$$ be the set of partitions of $t$ of $n$ parts and let $Y$ be the ${n \choose t}-$dimensional vector space over $\mathbb{F}_{p}$ with basis
$$\{c_{e_{0},e_{1},\ldots,e_{n}}\in \mathbb{F}_{q}: \, (e_{0},e_{1},\ldots, e_{n}) \in E^{t}_{n}\}.$$

Let's call $\mathcal{V}^{t}_{n}$ the Veronese image under the Veronese mapping given by:

$$\mathbb{F}_{p}(\sum_{i=0}^{n}x_{i}b_{i})\rightarrow \mathbb{F}_{p}(\sum_{E^{t}_{n}}c_{e_{0},\ldots,e_{n}}x^{e_{0}}x_{1}^{e_{1}}\cdots c^{e_{n}}_{n}),\ \ \ x_{i}\in \mathbb{F}_{p}.$$

The Veronese image of each $r-$dimensional subspace of $PG(n,p)$ is a sub-Veronesean variety $\mathcal{V}^{t}_{r}$ of $\mathcal{V}^{t}_{n}$, and all those subspaces are indexed by partitions in $\mathcal{P}(t)$. Thus by a Theorem due to Gmainer are invariant under the collineation group of the normal rational curve, (see \cite{Hav1}).

The $k-$rational points  $(p_{0},p_{1},\cdots, p_{n})$ of the normal rational curve $C$ correspond  to collinear points on $C$, 
that are defined over some Galois extension $l$ of $k$ and permuted by $\gal(l/k)$.

\cqd

\subsection{An application: three-point codes on the rational normal curve}


As we showed in Proposition \ref{codeNRC}, each subspace invariant under collineation of the NRC is indexed by a partition $\lambda \in \mathcal{P}(d)$. Let us call the base point associated to the partition $\lambda$ as $P_{\lambda}$.

\begin{theorem}\label{class}
Let $\sigma_{1}, \sigma_{2}, \sigma_{3}$ be three generators for the symmetric group $S_{d}$ and let $\lambda_{1}, \lambda_{2}$ and $\lambda_{3}$ be the partitions of $d$ indexing the corresponding irreducible representations in the special linear group $SL(n,\mathbb{F}_{q})$. Then any algebraic code defined over the NRC is covered by a divisor defined as linear combination 
of the base points $(P_{\lambda_{i}})_{1\leq i \leq 3}$ on the NRC, where the $\lambda_{i}$ are LR coefficients.
\end{theorem}

{\it Proof.}
Consider the divisors associated to the rational maps $f(x,y,z)=nx+my+lz$ defined over the rational normal curve $C_{d}$ defined over $\mathbb{F}_{q}$,  with $n, m$ and $l$ integer numbers. In particular, if $d | \, q^{2}-1$, the points $P=(\alpha, 0,0)$, $Q=(0,\beta,0)$ and $R=(0,0,\gamma)$ with $\alpha^{d}=1$, $\beta^{d}=1$ and $\gamma^{d}=1$, are $\mathbb{F}_{q^{2}}-$rational points on $C_{d}$, and the divisors $nP$, $mQ$ and $lR$ define codes on it.
Reciprocally, given a code on the NRC, by Proposition \ref{div}, the corresponding divisor defining the code is defined by a finite subgroup in the symmetric group. Since the symmetric group is generated by the 3 elements $\sigma_{1}, \sigma_{2}$ and $\sigma_{3}$, the divisor is a linear combination of the base points $(P_{\lambda_{i}})_{1\leq i \leq 3}$ on the NRC.

\cqd

\section{APPENDIX A: Explicit presentation of 3-point codes}
Given sets $I,J,K \subset \{0,1,\ldots, n\}$, of cardinality $r$, we can associate to them partitions $\lambda, \mu$ and $\nu$ as follows.
Let  $I=\{i_{1}<\ldots, < i_{r}\}\subset \{1,\ldots, n\}$ then the corresponding partition is defined as $\lambda=(i_{r}-r,\ldots, i_{1}-1)$. We consider the corresponding codes defined by the base points $c_{\lambda}$, $c_{\mu}$ and $c_{\nu}$, whenever the corresponding Littlewood-Richardson coefficient $c^{\nu}_{\lambda,\mu}$ is positive.
Next, we give an algorithm to compute the Littlewood-Richardson coefficients $c^{\nu}_{\lambda,\mu}$.
Horn defined sets of triples $(I,J,K)$ by the following inductive procedure (see \cite{Fu}):
$$U^{n}_{r}=\{(I,J,K)|\, \sum_{i\in I}+ \sum_{j\in J}=\sum_{k\in K}k+r(r+1)/2\},$$
$$T^{n}_{r}=\{(I,J,K)\in U^{n}_{r}|\, for\ all\ p<r \ and \ all \ (F,G,H) \in T^{r}_{p},$$
$$\sum_{f\in F}i_{f}+\sum_{g\in G}j_{g}\leq \sum_{h\in H}k_{h}+p(p+1)/2\}.$$

Note that Horn's algorithm  produces all the triples from the lowest values. Even if it is possible to start with a random generator set $I$, you need first to compute the lower values.
As a consequence of the classification Theorem \ref{class}, 
for any triple $(I,J,K)$ of indices admissible for the Horn problem the polynomials defined by $f(x)=\prod_{j\in J}(x-\alpha^{j}), g(x)=\prod_{i\in I}(x-\alpha^{i})$, and $ h(x)=\prod_{k\in K}(x-\alpha^{k})$ where $\alpha$ is a primitive element of $\mathbb{F}_{q^{m}}$ and $m$ is the least integer such that $n+1|\, p^{m}-1$ constitute a set of generators for the ideal of the corresponding algebraic code in the module of $n+1$ $\mathbb{F}_{q^{m}}$-rational points lying on the NRC.

Here we present a Sage \cite{SAGE} code calculating the $U^{n}_{r}$ and $T^{n}_{r}$ index sets, followed by a table containing all the cases till $n=4$ and $r=3$. The algorithm is implemented using Python: this involves calculate and iterate through $r-$combination of $n-$element. The running time is $O ({n \choose r}^{3} )$.

\newpage
{\tt
\begin{samepage}
\begin{verbatim}

from sage.combinat.subset import Subsets

def simple_cache(func):
    cache = dict()
    def cached_func(*args):
        if args not in cache:
            cache[args] = func(*args)
        return cache[args]
    cached_func.cache = cache
    return cached_func


@simple_cache
def getUnr(n, r):
    if r >= n:
        raise ValueError("r must be less than n: (n, r) = (%d, %d)" %(n, r))
    s = Subsets(range(1, n + 1), r)
    candidates = [(x, y, z) for x in s for y in s for z in s]
    return [tuple(map(sorted, (x, y, z))) for (x, y, z) in candidates if (
        sum(x) + sum(y)) == (sum(z) + r * (r + 1)/2)]


def index_filter(sub_index, index):
    if max(sub_index) > len(index):
        raise ValueError("%s must be valid indexes for %s" % (sub_index, index))
    # our indexes lists start at 1
    return [index[i - 1] for i in sub_index]

def condition((f, g, h), (i, j, k)):
    p = len(f)
    return sum(index_filter(f, i)) + sum(index_filter(g, j)) <= sum(
        index_filter(h, k)) + p*(p + 1)/2

def genTillR(r):
    return [getTnr(r, p) for p in range(1, r)]

@simple_cache
def getTnr(n, r):
    if r == 1:
        return getUnr(n, 1)
    else:
        return [(i, j, k) for (i, j, k) in getUnr(n, r) if all(
            all(condition((f, g, h), (i, j, k)) for (f, g, h) in triplets)
            for triplets in genTillR(r))]

\end{verbatim} 
\end{samepage}
} 
Here we list code's remarks
\begin{itemize}
\item the {\tt sorted()} mapping function in {\tt getUnr()} is necessary because the order of elements in {\tt Subsets} is unknown;
\item there is a 1-offset between index in Python lists and index sets we use;
\item the recursion in {\tt getTnr()} is factored out in {\tt getTillR()} call;
\item the cache decorator mitigate the perils of performing several times the same calculation in a function that already heavily recursive;
\item results are limiteted by constraints Python has on recursive function calls;
\item the fltering performed on $U^{n}_{r}$ to get $T^{n}_{r}$ is implemented by two nested calls to {\tt all()}
\end{itemize}

\begin{center}
\begin{tabular}{ | l | p{6cm} | p{6cm} | }

\hline
$(n, r)$   & $U^{n}_{r}$ &  $T^{n}_{r}$ \\
\hline
\hline
(2, 1) & $(\{1\}, \{1\}, \{1\})$, $(\{1\}, \{2\}, \{2\})$, $(\{2\}, \{1\}, \{2\})$ & $(\{1\}, \{1\}, \{1\})$, $(\{1\}, \{2\}, \{2\})$, $(\{2\}, \{1\}, \{2\})$ \\
\hline
\hline
(3, 1) & $(\{1\}, \{1\}, \{1\})$, $(\{1\}, \{2\}, \{2\})$, $(\{1\}, \{3\}, \{3\})$, $(\{2\}, \{1\}, \{2\})$, $(\{2\}, \{2\}, \{3\})$, $(\{3\}, \{1\}, \{3\})$ & $(\{1\}, \{1\}, \{1\})$, $(\{1\}, \{2\}, \{2\})$, $(\{1\}, \{3\}, \{3\})$, $(\{2\}, \{1\}, \{2\})$, $(\{2\}, \{2\}, \{3\})$, $(\{3\}, \{1\}, \{3\})$ \\
\hline
(3, 2) & $(\{1, 2\}, \{1, 2\}, \{1, 2\})$, $(\{1, 2\}, \{1, 3\}, \{1, 3\})$, $(\{1, 2\}, \{2, 3\}, \{2, 3\})$, $(\{1, 3\}, \{1, 2\}, \{1, 3\})$, $(\{1, 3\}, \{1, 3\}, \{2, 3\})$, $(\{2, 3\}, \{1, 2\}, \{2, 3\})$ & $(\{1, 2\}, \{1, 2\}, \{1, 2\})$, $(\{1, 2\}, \{1, 3\}, \{1, 3\})$, $(\{1, 2\}, \{2, 3\}, \{2, 3\})$, $(\{1, 3\}, \{1, 2\}, \{1, 3\})$, $(\{1, 3\}, \{1, 3\}, \{2, 3\})$, $(\{2, 3\}, \{1, 2\}, \{2, 3\})$ \\
\hline
\hline
(4, 1)  & $(\{1\}, \{1\}, \{1\})$, $(\{1\}, \{2\}, \{2\})$, $(\{1\}, \{3\}, \{3\})$, $(\{1\}, \{4\}, \{4\})$, $(\{2\}, \{1\}, \{2\})$, $(\{2\}, \{2\}, \{3\})$, $(\{2\}, \{3\}, \{4\})$, $(\{3\}, \{1\}, \{3\})$, $(\{3\}, \{2\}, \{4\})$, $(\{4\}, \{1\}, \{4\})$ & $(\{1\}, \{1\}, \{1\})$, $(\{1\}, \{2\}, \{2\})$, $(\{1\}, \{3\}, \{3\})$, $(\{1\}, \{4\}, \{4\})$, $(\{2\}, \{1\}, \{2\})$, $(\{2\}, \{2\}, \{3\})$, $(\{2\}, \{3\}, \{4\})$, $(\{3\}, \{1\}, \{3\})$, $(\{3\}, \{2\}, \{4\})$, $(\{4\}, \{1\}, \{4\})$ \\
\hline
(4, 2) & $(\{1, 2\}, \{1, 2\}, \{1, 2\})$, $(\{1, 2\}, \{1, 3\}, \{1, 3\})$, $(\{1, 2\}, \{1, 4\}, \{1, 4\})$, $(\{1, 2\}, \{1, 4\}, \{2, 3\})$, $(\{1, 2\}, \{2, 3\}, \{1, 4\})$, $(\{1, 2\}, \{2, 3\}, \{2, 3\})$, $(\{1, 2\}, \{2, 4\}, \{2, 4\})$, $(\{1, 2\}, \{3, 4\}, \{3, 4\})$, $(\{1, 3\}, \{1, 2\}, \{1, 3\})$, $(\{1, 3\}, \{1, 3\}, \{1, 4\})$, $(\{1, 3\}, \{1, 3\}, \{2, 3\})$, $(\{1, 3\}, \{1, 4\}, \{2, 4\})$, $(\{1, 3\}, \{2, 3\}, \{2, 4\})$, $(\{1, 3\}, \{2, 4\}, \{3, 4\})$, $(\{1, 4\}, \{1, 2\}, \{1, 4\})$, $(\{1, 4\}, \{1, 2\}, \{2, 3\})$, $(\{1, 4\}, \{1, 3\}, \{2, 4\})$, $(\{1, 4\}, \{1, 4\}, \{3, 4\})$, $(\{1, 4\}, \{2, 3\}, \{3, 4\})$, $(\{2, 3\}, \{1, 2\}, \{1, 4\})$, $(\{2, 3\}, \{1, 2\}, \{2, 3\})$, $(\{2, 3\}, \{1, 3\}, \{2, 4\})$, $(\{2, 3\}, \{1, 4\}, \{3, 4\})$, $(\{2, 3\}, \{2, 3\}, \{3, 4\})$, $(\{2, 4\}, \{1, 2\}, \{2, 4\})$, $(\{2, 4\}, \{1, 3\}, \{3, 4\})$, $(\{3, 4\}, \{1, 2\}, \{3, 4\})$ & $(\{1, 2\}, \{1, 2\}, \{1, 2\})$, $(\{1, 2\}, \{1, 3\}, \{1, 3\})$, $(\{1, 2\}, \{1, 4\}, \{1, 4\})$, $(\{1, 2\}, \{2, 3\}, \{2, 3\})$, $(\{1, 2\}, \{2, 4\}, \{2, 4\})$, $(\{1, 2\}, \{3, 4\}, \{3, 4\})$, $(\{1, 3\}, \{1, 2\}, \{1, 3\})$, $(\{1, 3\}, \{1, 3\}, \{1, 4\})$, $(\{1, 3\}, \{1, 3\}, \{2, 3\})$, $(\{1, 3\}, \{1, 4\}, \{2, 4\})$, $(\{1, 3\}, \{2, 3\}, \{2, 4\})$, $(\{1, 3\}, \{2, 4\}, \{3, 4\})$, $(\{1, 4\}, \{1, 2\}, \{1, 4\})$, $(\{1, 4\}, \{1, 3\}, \{2, 4\})$, $(\{1, 4\}, \{1, 4\}, \{3, 4\})$, $(\{2, 3\}, \{1, 2\}, \{2, 3\})$, $(\{2, 3\}, \{1, 3\}, \{2, 4\})$, $(\{2, 3\}, \{2, 3\}, \{3, 4\})$, $(\{2, 4\}, \{1, 2\}, \{2, 4\})$, $(\{2, 4\}, \{1, 3\}, \{3, 4\})$, $(\{3, 4\}, \{1, 2\}, \{3, 4\})$ \\
\hline
(4, 3) & $(\{1, 2, 3\}, \{1, 2, 3\}, \{1, 2, 3\})$, $(\{1, 2, 3\}, \{1, 2, 4\}, \{1, 2, 4\})$, $(\{1, 2, 3\}, \{1, 3, 4\}, \{1, 3, 4\})$, $(\{1, 2, 3\}, \{2, 3, 4\}, \{2, 3, 4\})$, $(\{1, 2, 4\}, \{1, 2, 3\}, \{1, 2, 4\})$, $(\{1, 2, 4\}, \{1, 2, 4\}, \{1, 3, 4\})$, $(\{1, 2, 4\}, \{1, 3, 4\}, \{2, 3, 4\})$, $(\{1, 3, 4\}, \{1, 2, 3\}, \{1, 3, 4\})$, $(\{1, 3, 4\}, \{1, 2, 4\}, \{2, 3, 4\})$, $(\{2, 3, 4\}, \{1, 2, 3\}, \{2, 3, 4\})$  & $(\{1, 2, 3\}, \{1, 2, 3\}, \{1, 2, 3\})$, $(\{1, 2, 3\}, \{1, 2, 4\}, \{1, 2, 4\})$, $(\{1, 2, 3\}, \{1, 3, 4\}, \{1, 3, 4\})$, $(\{1, 2, 3\}, \{2, 3, 4\}, \{2, 3, 4\})$, $(\{1, 2, 4\}, \{1, 2, 3\}, \{1, 2, 4\})$, $(\{1, 2, 4\}, \{1, 2, 4\}, \{1, 3, 4\})$, $(\{1, 2, 4\}, \{1, 3, 4\}, \{2, 3, 4\})$, $(\{1, 3, 4\}, \{1, 2, 3\}, \{1, 3, 4\})$, $(\{1, 3, 4\}, \{1, 2, 4\}, \{2, 3, 4\})$, $(\{2, 3, 4\}, \{1, 2, 3\}, \{2, 3, 4\})$ \\
\hline
\end{tabular}

\end{center}

\end{document}